 \documentclass[final,3p,times]{elsarticle}

\usepackage{amssymb}
 \usepackage{amsthm}
 \usepackage[T1]{fontenc}
 
\usepackage[ansinew]{inputenc} 
\usepackage{mathrsfs}
\usepackage{euscript}
\usepackage{natbib}

\newtheorem{theorem}{Theorem}
\newtheorem{lemma}{Lemma}
\newtheorem{proposition}{Proposition}
\newdefinition{rmk}{ Remark}
\newproof{pf}{Proof of Theorem 1}

\begin{document}

\begin{frontmatter}



\title{ Finite time blow up for wave equations with strong damping in an exterior domain}

 \author[label]{Ahmad Z. FINO}
 \address[label]{LaMA-Liban, Lebanese University, Faculty of Sciences, Department of Mathematics, P.O. Box 826 Tripoli, Lebanon}


 \ead{ahmad.fino01@gmail.com; afino@ul.edu.lb}
 \begin{abstract}
We consider the initial boundary value problem in exterior domain for strongly damped wave equations with power-type nonlinearity $|u|^p$. We will establish blow-up results under some conditions on the initial data and the exponent $p$.
\end{abstract}


\end{frontmatter}


\section{Introduction}
\setcounter{equation}{0}

This paper concerns the initial boundary value problem of the strongly damped wave equation in an exterior domain. Let $\Omega\subset\mathbb{R}^n$ be an exterior domain whose obstacle 
$\mathcal{O}\subset\mathbb{R}^n$ is bounded with smooth compact boundary $\partial\Omega$. We consider the initial boundary value problem
\begin{equation}\label{eq1}
\left\{\begin{array}{ll}
\,\, \displaystyle {u_{tt}-\Delta
u -\Delta u_t =|u|^p} &\displaystyle {t>0,x\in {\Omega},}\\
{}\\
\displaystyle{u(0,x)=  u_0(x),\;\;u_t(0,x)=  u_1(x)\qquad\qquad}&\displaystyle{x\in {\Omega},}\\
{}\\
\displaystyle{u=0,\qquad\qquad}&\displaystyle{t\geq0,\;x\in {\partial\Omega},}
\end{array}
\right.
\end{equation} 
where the unknown function $u$ is real-valued, $n\geq 1$, and $p>1.$ Throughout this paper, we assume that 
\begin{equation}\label{condition1}
(u_0,u_1)\in (H^2(\Omega)\cap H^1_0(\Omega))\times L^2(\Omega),\quad\mbox{and}\quad u_0,u_1\geq0.
\end{equation}
Without loss of generality, we assume that $0\in\mathcal{O}\subset\subset B(R)$, where $B(R):=\{x\in\mathbb{R}^n:\;|x|< R\}$ is a ball of radius $R$ centred at the origin.\\
For the simplicity of notations, $\|\cdotp\|_q$ and $\|\cdotp\|_{H^1}$ $(1\leq q\leq \infty)$ stand for the usual $L^q(\Omega)$-norm and $H^1_0(\Omega)$-norm, respectively.

First, the following local well-posedness result is needed. 
\begin{proposition}\label{prop1}\cite[see Proposition 2.1]{IkehataInoue}\\
Let $1<p<\infty$ for $n=1,2$ and $1<p\leq \frac{n}{n-2}$ for $n\geq3$. Under the assumption $(\ref{condition1})$, there exists a maximal existence time $T_{\max}>0$ such that the problem $(\ref{eq1})$ possesses a unique weak solution
$$u\in C([0,T_{\max}),H^1_0(\Omega))\cap C^1([0,T_{\max}), L^2(\Omega)),$$
where $0< T_{\max}\leq\infty.$ In addition:
\begin{equation}\label{alternative}
\mbox{either}\;\;T_{\max}=\infty \quad\mbox{or else}\quad T_{\max}<\infty \;\;\mbox{and}\;\; \|u(t,\cdotp)\|_{H^1}+\|u_t(t,\cdotp)\|_2\rightarrow\infty\;\;\mbox{as}\;\; t\rightarrow T_{\max}.
 \end{equation}
\end{proposition}
\begin{rmk}
We say that $u$ is a global solution of $(\ref{eq1})$ if $T_{\max}=\infty,$ while in the case of $T_{\max}<\infty,$ we say that $u$ blows up in finite time.
\end{rmk}
Our main result is the following
\begin{theorem}\label{theo1}
 Assume that the initial data satisfies $(\ref{condition1})$ such that
 $$\int_{\Omega}[u_1(x)-\Delta u_0(x)]\phi_0(x)\,dx>0,$$
 where $\phi_0(x)$ is defined below (see Lemma \ref{lemma1}, \ref{lemma4}, \ref{lemma2}). If 
 $$\left\{\begin{array}{l}
\displaystyle 1<p\leq 3,\quad \mbox{for}\; n=1,\\
\displaystyle   1<p< 3,\quad \mbox{for}\;n=2,\\
  \displaystyle  1<p\leq 1+\frac{2}{n-1},\quad \mbox{for}\;n\geq3,
 \end{array}
 \right.
 $$
then the solution of the problem $(\ref{eq1})$ blows up in finite time.
\end{theorem}

This paper is organized as follows: in Section \ref{sec2}, we present several preliminaries. Section \ref{sec3} contains the proofs of the blow-up theorem (Theorem \ref{theo1}).


\section{Preliminaries}\label{sec2}
\setcounter{equation}{0}
In this section, we give some preliminary properties that will be used in the proof of Theorem \ref{theo1}.
\begin{lemma}\label{lemma1}
There exists a function $\phi_0(x)\in C^2(\Omega)\cap C(\overline{\Omega})$ satisfying the following boundary value problem
\begin{equation}\label{function1}
\left\{\begin{array}{l}
\Delta \phi_0(x)=0, \;\mbox{in}\;\Omega,\;\;n\geq 3,\\
\phi_0|_{\partial\Omega}=0,\\
|x|\rightarrow\infty,\quad \phi_0(x)\rightarrow 1.
\end{array}
\right.
\end{equation}
Moreover, $\phi_0(x)$ satisfies: 
\begin{itemize}
 \item $0<\phi_0(x)<1$, for all $x\in\Omega$.
 \item $\phi_0(x)\geq C$, for all $|x|\gg1$.
 \item $|\nabla \phi_0(x)|\leq \frac{C}{|x|^{n-1}}$, for all $|x|\gg1$.
 \end{itemize}
\end{lemma}
\proof From \cite[Lemma 2.2]{ZhouHan} there exists a regular solution $\phi_0$ of (\ref{function1}) such that $0<\phi_0(x)<1$, for all $x\in\Omega$. To obtain the last two properties of $\phi_0$, it is easy to see that since $\mathcal{O}$ is bounded, there exist $r_2>r_1>0$ such that $B_{r_1}\subseteq\mathcal{O}\subseteq B_{r_2}$, where $B_r$ stands for the open ball with center zero and radius $r$. By the maximum principle we conclude that $\phi_1(x)\leq \phi_0(x)\leq \phi_2(x)$ in $\Omega$, where $\phi_1(x)$ and $\phi_2(x)$ are, respectively, the solution of (\ref{function1}) on $\mathbb{R}^n\setminus B_{r_1}$ and $\mathbb{R}^n\setminus B_{r_2}$. We remember tha $\phi_i(x)=r_i^{2-n}-|x|^{2-n}$, $i=1,2$. Moreover, the standard elliptic theory implies that $|\nabla \phi_0(x)|\sim |\nabla \phi_{i}(x)|$, $i=1,2$. As $\phi_{1}(x)|\geq C$ and $|\nabla\phi_{i}(x)|\leq \frac{C}{|x|^{n-1}}$ when $|x|\gg1$, this complete the proof.$\hfill\square$\\

Similarly, we have the following
\begin{lemma}\label{lemma4}\cite[Lemma 2.5]{Han}
There exists a function $\phi_0(x)\in C^2(\Omega)\cap C(\overline{\Omega})$ satisfying the following boundary value problem
\begin{equation}\label{function2}
\left\{\begin{array}{l}
\Delta \phi_0(x)=0, \;\mbox{in}\;\Omega,\;\;n=2,\\
\phi_0|_{\partial\Omega}=0,\\
|x|\rightarrow\infty,\quad \phi_0(x)\rightarrow +\infty,\;\;\mbox{and}\;\phi_0(x)\;\mbox{increase at the rate of}\;\ln (|x|).
\end{array}
\right.
\end{equation}
Moreover, $\phi_0(x)$ satisfies: 
\begin{itemize}
 \item $0<\phi_0(x)\leq C\ln (|x|)$, for all $x\in\Omega$.
 \item $\phi_0(x)\geq C$, for all $|x|\gg1$.
 \item $|\nabla \phi_0(x)|\leq \frac{C}{|x|}$, for all $|x|\gg1$.
 \end{itemize}
\end{lemma}
\begin{lemma}\label{lemma2}\cite[Lemma 2.2]{Han2}
There exists a function $\phi_0(x)\in C^2([0,\infty))$ satisfying the following boundary value problem
\begin{equation}\label{function3}
\left\{\begin{array}{l}
\Delta \phi_0(x)=0, \;x>0,\\
\phi_0|_{x=0}=0,\\
x\rightarrow\infty,\quad \phi_0(x)\rightarrow +\infty,\;\;\mbox{and}\;\phi_0(x)\;\mbox{increase at the rate of linear function $x$}.
\end{array}
\right.
\end{equation}
Moreover, $\phi_0(x)$ satisfies: there exist two positive constants $C_1$ and $C_2$ such that, for all $x>0$, we have $C_1x\leq \phi_0(x)\leq C_2x$. In fact, we can take $\phi_0(x)=Cx$.
\end{lemma}


\section{Proof of Theorem \ref{theo1}}\label{sec3}
\setcounter{equation}{0}

\pf
The idea to prove Theorem \ref{theo1} is to use the variational formulation of the weak solution by choosing the appropriate test function. Note that the harmonic functions in Lemma \ref{lemma1}, \ref{lemma4} and \ref{lemma2} play a crucial role in the exterior domain, because of their good behaviour and vanishing on the boundary $\partial\Omega$.

We argue by contradiction assuming that $u$ is 
not a blow-up solution of (\ref{eq1}), we have
 \begin{eqnarray}\label{newweaksolution}
&{}&\int_0^T\int_{\Omega}|u|^{p}\varphi\,dx\,dt+\int_{\Omega}[u_1(x)-\Delta u_0(x)]\varphi(0,x)\,dx-\int_{\Omega}u_0(x)\varphi_t(0,x)\,dx\nonumber\\
&{}&=\int_0^T\int_{\Omega}u\varphi_{tt}\,dx\,dt+\int_0^T\int_{\Omega}u\Delta\varphi_t\,dx\,dt-\int_0^T\int_{\Omega}u\Delta\varphi\,dx\,dt,
\end{eqnarray}
for all $T>0$ and all compactly supported function $\varphi\in C^2([0,T]\times\Omega)$ 
such that $\varphi(\cdotp,T)=0$ and $\varphi_t(\cdotp,T)=0$. Take $\varphi(x,t)=\phi_0(x)\varphi^\ell_T(x)\eta_T^{k}(t)$ where $\phi_0$ is the harmonic function introduced in Lemma \ref{lemma1}, \ref{lemma4} and \ref{lemma2}, $\eta_T(t):=\eta(\frac{t^2}{T^2})$, $\ell, k\gg1$, and $\eta(\cdotp)\in C^\infty(\mathbb{R}_+)$ is a cut-off non-increasing function such that
$$\eta(t):=\left\{\begin {array}{ll}\displaystyle{1}&\displaystyle{\quad\mbox{if }0\leq t\leq 1/4}\\
\displaystyle{0}&\displaystyle{\quad\mbox {if }t\geq 1,}
\end {array}\right.$$
$0\leq \eta(t) \leq 1$ and $|\eta^{'}(t)|\leq C$ for some $C>0$ and all $t>0$; and $\varphi_T(x)=\Phi(\frac{|x|}{T})$ with the following smooth, non-increasing cut-off function
$$\Phi(r):=\left\{\begin {array}{ll}\displaystyle{1}&\displaystyle{\quad\mbox{if }0\leq r\leq 1}\\
\displaystyle{0}&\displaystyle{\quad\mbox {if }r\geq 2,}
\end {array}\right.$$\\
such that $0\leq\Phi(r)\leq 1$, $|\Phi'(r)|\leq C/r$ and $|\Phi''(r)|\leq C/r^2$. We obtain
 \begin{eqnarray}\label{weaksolution2}
&{}&\int_0^T\int_{\Omega_1}|u|^{p}\varphi\,dx\,dt+\int_{\Omega_1}[u_1(x)-\Delta u_0(x)]\phi_0(x)\varphi_T^\ell(x)\,dx\nonumber\\
&{}&=\int_0^T\int_{\Omega_1}u\phi_0(x)\varphi_T^\ell(x)\partial_t^2(\eta^k_T(t))\,dx\,dt+\int_0^T\int_{\Omega_1}u\Delta[\phi_0(x)\varphi_T^\ell(x)]\partial_t(\eta^k_T(t))\,dx\,dt-\int_0^T\int_{\Omega_1}u\Delta[\phi_0(x)\varphi_T^\ell(x)]\eta^k_T(t)\,dx\,dt\nonumber\\
&{}&=:I_1+I_2+I_3
\end{eqnarray}
where $\Omega_1:=\{x\in\Omega;\;|x|\leq 2T\}$. At this stage, we have to distinguishes three cases:\\

\noindent $\bullet$ \underline{Case 1: $n\geq 3$}. To estimate the right-hand side of (\ref{weaksolution2}), we introduce the term 
$\varphi^{1/p}\varphi^{-1/p}$ in $I_1$, and we use Young's inequality to obtain
\begin{eqnarray}\label{I1}
I_1&\leq&\int_0^T\int_{\Omega_1}|u|\;\varphi^{1/p}\varphi^{-1/p}\phi_0(x)\varphi^\ell_B(x)\left|\partial_t^2[\eta_T^{k}(t)]\right|\,dx\,dt\nonumber\\
&\leq&\frac{1}{6} \int_0^T\int_{\Omega_1}|u|^p\;\varphi \,dx\,dt+C\int_0^T\int_{\Omega_1}\varphi^{-p'/p}\phi_0^{p'}(x)\varphi^{\ell p'}_T(x)\left|\partial_t^2[\eta_T^{k}(t)]\right|^{p'}\,dx\,dt\nonumber\\
&\leq&\frac{1}{6} \int_0^T\int_{\Omega_1}|u|^p\;\varphi \,dx\,dt+C\int_0^T\int_{\Omega_1}\phi_0(x)\varphi_T^\ell(x)\eta_T(t)^{(k-2)p'}\left|\partial_t \eta_T(t)\right|^{2p'}\,dx\,dt\nonumber\\&{}&+C\int_0^T\int_{\Omega_1}\phi_0(x)\varphi_T^\ell(x)\eta_T(t)^{(k-1)p'}\left|\partial^2_t \eta_T(t)\right|^{p'}\,dx\,dt.
\end{eqnarray}
On the other hand, using Lemma \ref{lemma1} with all properties of $\phi_0$, $T\gg1$, and Young's inequality, we conclude that
\begin{eqnarray}\label{I2}
I_2&\leq&C\int_0^T\int_{\Omega_1}|u|\varphi_T^{\ell-1}(x)
\left|\nabla\phi_0(x)\right|\left|\nabla\varphi_T(x)\right||\partial_t(\eta^k_T(t))|\,dx\,dt\nonumber\\
&{}&+\,C\int_0^T\int_{\Omega_1}|u|\varphi_T^{\ell-2}(x)
\phi_0(x)\left|\nabla\varphi_T(x)\right|^2|\partial_t(\eta^k_T(t))|\,dx\,dt\nonumber\\
&{}&+\,C\int_0^T\int_{\Omega_1}|u|\varphi_T^{\ell-1}(x)\phi_0(x)
\left|\Delta\varphi_T(x)\right||\partial_t(\eta^k_T(t))|\,dx\,dt\nonumber\\
&=&C\int_0^T\int_{\Omega_1}|u|\;\varphi^{1/p}\varphi^{-1/p}\varphi_T^{\ell-1}(x)
\left|\nabla\phi_0(x)\right|\left|\nabla\varphi_T(x)\right||\partial_t(\eta^k_T(t))|\,dx\,dt\nonumber\\
&{}&+\,C\int_0^T\int_{\Omega_1}|u|\;\varphi^{1/p}\varphi^{-1/p}\varphi_T^{\ell-2}(x)
\phi_0(x)\left|\nabla\varphi_T(x)\right|^2|\partial_t(\eta^k_T(t))|\,dx\,dt\nonumber\\
&{}&+\,C\int_0^T\int_{\Omega_1}|u|\;\varphi^{1/p}\varphi^{-1/p}\varphi_T^{\ell-1}(x)\phi_0(x)
\left|\Delta\varphi_T(x)\right||\partial_t(\eta^k_T(t))|\,dx\,dt\nonumber\\
&\leq&\frac{1}{6} \int_0^T\int_{\Omega_1}|u|^p\;\varphi \,dx\,dt+C\int_0^T\int_{\nabla\Omega_1}\varphi_T^{\ell-p'}(x)\eta_T^{k-p'}(t)|\nabla\phi_0(x)|^{p'}\left|\nabla\varphi_T(x)\right|^{p'}|\partial_t(\eta_T(t))|^{p'}\,dx\,dt\nonumber\\
&{}&+\,C\int_0^T\int_{\nabla\Omega_1}\varphi_T^{\ell-2p'}(x)\eta_T^{k-p'}(t)
\left|\nabla\varphi_T(x)\right|^{2p'}|\partial_t(\eta_T(t))|^{p'}\,dx\,dt\nonumber\\
&{}&+\,C\int_0^T\int_{\nabla\Omega_1}\varphi_T^{\ell-p'}(x)\eta_T^{k-p'}(t)
\left|\Delta\varphi_T(x)\right|^{p'}|\partial_t(\eta_T(t))|^{p'}\,dx\,dt,
\end{eqnarray}
where $\nabla\Omega_1:=\{x\in\Omega;\;T\leq |x|\leq 2T\}$. Similarly,
\begin{eqnarray}\label{I3}
I_3&\leq&\frac{1}{6} \int_0^T\int_{\Omega_1}|u|^p\;\varphi \,dx\,dt+C\int_0^T\int_{\nabla\Omega_1}\varphi_T^{\ell-p'}(x)\eta_T^{k}(t)|\nabla\phi_0(x)|^{p'}\left|\nabla\varphi_T(x)\right|^{p'}\,dx\,dt\nonumber\\
&{}&+\,C\int_0^T\int_{\nabla\Omega_1}\varphi_T^{\ell-2p'}(x)\eta_T^{k}(t)
\left|\nabla\varphi_T(x)\right|^{2p'}\,dx\,dt\nonumber\\
&{}&+\,C\int_0^T\int_{\nabla\Omega_1}\varphi_T^{\ell-p'}(x)\eta_T^{k}(t)
\left|\Delta\varphi_T(x)\right|^{p'}\,dx\,dt,
\end{eqnarray}
Using (\ref{I1})-(\ref{I3}), it follows from (\ref{weaksolution2}) that
\begin{eqnarray}\label{weaksolution3}
&{}&\int_{\Omega_1}[u_1(x)-\Delta u_0(x)]\phi_0(x)\varphi_T^\ell(x)\,dx\nonumber\\
&{}& \leq \frac{1}{2}\int_0^T\int_{\Omega_1}|u|^p\;\varphi \,dx\,dt+\int_{\Omega_1}[u_1(x)-\Delta u_0(x)]\phi_0(x)\varphi_T^\ell(x)\,dx\nonumber\\
&{}& \leq C\int_0^T\int_{\Omega_1}\phi_0(x)\varphi_T^\ell(x)\eta_T(t)^{(k-2)p'}\left|\partial_t \eta_T(t)\right|^{2p'}\,dx\,dt\nonumber\\
&{}&+\,C\int_0^T\int_{\Omega_1}\phi_0(x)\varphi_T^\ell(x)\eta_T(t)^{(k-1)p'}\left|\partial^2_t \eta_T(t)\right|^{p'}\,dx\,dt\nonumber\\
&{}&+\;C\int_0^T\int_{\nabla\Omega_1}\varphi_T^{\ell-p'}(x)\eta_T^{k-p'}(t)|\nabla\phi_0(x)|^{p'}\left|\nabla\varphi_T(x)\right|^{p'}|\partial_t(\eta_T(t))|^{p'}\,dx\,dt\nonumber\\
&{}&+\,C\int_0^T\int_{\nabla\Omega_1}\varphi_T^{\ell-2p'}(x)\eta_T^{k-p'}(t)
\left|\nabla\varphi_T(x)\right|^{2p'}|\partial_t(\eta_T(t))|^{p'}\,dx\,dt\nonumber\\
&{}&+\,C\int_0^T\int_{\nabla\Omega_1}\varphi_T^{\ell-p'}(x)\eta_T^{k-p'}(t)
\left|\Delta\varphi_T(x)\right|^{p'}|\partial_t(\eta_T(t))|^{p'}\,dx\,dt\nonumber\\
&{}&+\,C\int_0^T\int_{\nabla\Omega_1}\varphi_T^{\ell-p'}(x)\eta_T^{k}(t)|\nabla\phi_0(x)|^{p'}\left|\nabla\varphi_T(x)\right|^{p'}\,dx\,dt\nonumber\\
&{}&+\,C\int_0^T\int_{\nabla\Omega_1}\varphi_T^{\ell-2p'}(x)\eta_T^{k}(t)
\left|\nabla\varphi_T(x)\right|^{2p'}\,dx\,dt+\,C\int_0^T\int_{\nabla\Omega_1}\varphi_T^{\ell-p'}(x)\eta_T^{k}(t)
\left|\Delta\varphi_T(x)\right|^{p'}\,dx\,dt,
\end{eqnarray}
Now, we have to distinguishes 2 subcases.\\ 
\noindent $\bullet$ \underline{Case (i): $1<p<1+\frac{2}{n-1}$}. By Lemma \ref{lemma1}, we have: $|\nabla\phi_0(x)|\leq\frac{C}{|x|^{n-1}}\leq \frac{C}{T^{n-1}}\leq \frac{C}{T}$) in $\nabla\Omega_1$, therefore, using the change of variables: $y=T^{-1}x,$ $s=T^{-1}t$, we get from $(\ref{weaksolution3})$ that 
\begin{eqnarray}\label{weaksolution4}
\int_{\Omega_1}[u_1(x)-\Delta u_0(x)]\phi_0(x)\varphi_T^\ell(x)\,dx&\leq&C\;T^{-2p'+1+n}+C\;T^{-3p'+1+1n}\nonumber\\
&\leq&C\;T^{-2p'+1+n},
\end{eqnarray}
where $C$ is independent of $T$. As $p<1+\frac{2}{n-1}\Longleftrightarrow -2p'+1+n<0$, it follows, by letting $T\rightarrow\infty$ that
$$0<\int_{\Omega}[u_1(x)-\Delta u_0(x)]\phi_0(x)\,dx\leq0;$$
contradiction.\\
\noindent $\bullet$ \underline{Case (ii): $p=1+\frac{2}{n-1}$}. From (\ref{weaksolution3}) in the Case 1 and the fact that $p=1+\frac{2}{n-1},$ there exists a positive constant $D$ independent of $T$ such that
$$\int_0^{T}\int_{\Omega_1}|u|^p\;\varphi \,dx\,dt\leq D,\qquad\mbox{for all}\; T>0,$$
which implies that
\begin{equation}\label{lefthand}
\int_{T/2}^{T}\int_{\Omega_1}|u|^p\;\varphi \,dx\,dt,\;\int_{T/2}^{T}\int_{\nabla\Omega_1}|u|^p\;\varphi \,dx\,dt,\;\int_{0}^{T}\int_{\nabla\Omega_1}|u|^p\;\varphi \,dx\,dt\rightarrow0\quad\mbox{as}\;\;T\rightarrow\infty.
\end{equation}
On the other hand, we use H\"older's inequality instead of Young's one in $I_1$, $I_2$, and $I_3$, together with the same change of variables, we get
 \begin{eqnarray}\label{newI1}
I_1&\leq&\left(\int_{T/2}^T\int_{\Omega_1}|u|^p\;\varphi \,dx\,dt\right)^{1/p}\left(C\int_0^T\int_{\Omega_1}\phi_0(x)\varphi_T^\ell(x)\left[\eta_T(t)^{(k-2)p'}\left|\partial_t \eta_T(t)\right|^{2p'}+\eta_T(t)^{(k-1)p'}\left|\partial^2_t \eta_T(t)\right|^{p'}\right]\,dx\,dt\right)^{1/p'}\nonumber\\
&\leq&C\;T^{-2+\frac{1+n}{p'}}\left(\int_{T/2}^T\int_{\Omega_1}|u|^p\;\varphi \,dx\,dt\right)^{1/p}\nonumber\\
&=&C\;\left(\int_{T/2}^T\int_{\Omega_1}|u|^p\;\varphi \,dx\,dt\right)^{1/p},
\end{eqnarray}
thanks to the fact that $p=1+\frac{2}{n-1}$. Similarly 
\begin{eqnarray}\label{newI2}
I_2&\leq&C\;\left(\int_{T/2}^T\int_{\nabla\Omega_1}|u|^p\;\varphi \,dx\,dt\right)^{1/p},
\end{eqnarray}
and
\begin{eqnarray}\label{newI3}
I_3&\leq&C\;\left(\int_{0}^T\int_{\nabla\Omega_1}|u|^p\;\varphi \,dx\,dt\right)^{1/p}.
\end{eqnarray}
Finally, using (\ref{newI1})-(\ref{newI3}), it follows from (\ref{weaksolution2}) that
\begin{eqnarray*}
\int_{\Omega_1}[u_1(x)-\Delta u_0(x)]\phi_0(x)\varphi_T^\ell(x)\,dx&\leq& C\;\left(\int_{T/2}^T\int_{\Omega_1}|u|^p\;\varphi \,dx\,dt\right)^{1/p}\\
&{}&+C\;\left(\int_{T/2}^T\int_{\nabla\Omega_1}|u|^p\;\varphi \,dx\,dt\right)^{1/p}+C\;\left(\int_{0}^T\int_{\nabla\Omega_1}|u|^p\;\varphi \,dx\,dt\right)^{1/p},
\end{eqnarray*}
hence, by letting $T\rightarrow\infty$ and using (\ref{lefthand}), we get a contradiction.\\

\noindent $\bullet$ \underline{Case 2: $n=2$}. In this case, we have a blow-up result just in the sub-critical case ($1<p<1+\frac{2}{n-1}=3$). By repeating the same calculation in the Case of $n\geq3$ and using Lemma \ref{lemma4} instead of Lemma \ref{lemma1} (noted that the big difference is the fact that $\phi_0(x)\leq C\ln (|x|)$), we easily conclude that
$$I_1\leq \frac{1}{6} \int_0^T\int_{\Omega_1}|u|^p\;\varphi \,dx\,dt+ C\ln(T)\,T^{-2p'+3},
$$
$$I_2\leq \frac{1}{6} \int_0^T\int_{\Omega_1}|u|^p\;\varphi \,dx\,dt+ C\,T^{-3p'+3}+C\,\ln(T)\,T^{-3p'+3},
$$
and
$$I_3\leq \frac{1}{6} \int_0^T\int_{\Omega_1}|u|^p\;\varphi \,dx\,dt+ C\,T^{-2p'+3}+C\,\ln(T)\,T^{-2p'+3}.$$
This implies that
$$\int_{\Omega_1}[u_1(x)-\Delta u_0(x)]\phi_0(x)\varphi_T^\ell(x)\,dx\leq C\ln(T)\,T^{-2p'+3}\leq C\,T^{-p'+3/2},$$
where we have used, e.g., the fact that $\ln(T)\leq C\,T^{p'-3/2}$. By letting $T$ goes to infinity and using $p<{3}$, we obtain the desired contradiction.\\

\noindent $\bullet$ \underline{Case 3: $n=1$}. For the case $1<p<3$, repeat the same calculation as in the Case of $n\geq3$ and using Lemma \ref{lemma2} instead of Lemma \ref{lemma1}, we easily get
$$I_1\leq \frac{1}{6} \int_0^T\int_{\Omega_1}|u|^p\;\varphi \,dx\,dt+ C\,T^{-2p'+3},
$$
$$I_2\leq \frac{1}{6} \int_0^T\int_{\Omega_1}|u|^p\;\varphi \,dx\,dt+C\,T^{-3p'+3},
$$
and
$$I_3\leq \frac{1}{6} \int_0^T\int_{\Omega_1}|u|^p\;\varphi \,dx\,dt+ C\,T^{-2 p'+3}.$$
Using the change of variables: $y=T^{-\alpha}x,$ $s=T^{-1}t$, we get from $(\ref{weaksolution3})$ that 
\begin{eqnarray*}
\int_{\Omega_1}[u_1(x)-\Delta u_0(x)]\phi_0(x)\varphi_T^\ell(x)\,dx\leq C\;T^{-2 p'+3},
\end{eqnarray*}
which leads to a contradiction by letting $T\rightarrow\infty$.\\
For the critical case $p=3$, we get the contradiction by applying a similar calculation as in the case (ii) above by taking into account the support of $\nabla\varphi_T$, $\Delta\varphi_T$, and $\partial_t\eta_T$.\\
This completes the proof of Theorem $\ref{theo1}$. $\hfill\square$




\section*{Acknowledgements}
The author would like to express sincere gratitude to Professor Ryo Ikehata for valuable discussion. 

\bibliographystyle{elsarticle-num}



\section*{References}
 \bibitem{Han} {Wei Han}, Concerning the Strauss Conjecture for the subcritical and
critical cases on the exterior domain in two space dimensions, Nonlinear Analysis ${\bf 84}$ $(2013),$ $136-145.$

\bibitem{Han2} {Wei Han}, Blow Up of Solutions to One Dimensional
Initial-Boundary Value Problems for Semilinear Wave Equations with Variable Coefficients, J. Part. Diff. Eq. ${\bf 26}$ $(2013),$ No. 2 $138-150.$ 

 \bibitem{IkehataInoue} {R. Ikehata, Yu-ki Inoue}, Global existence of weak solutions for two-dimensional semilinear wave equations with strong damping in an exterior domain, Nonlinear Analysis ${\bf 68}$ $(2008),$ $154-169$.

\bibitem{ZhouHan} {Yi Zhou, Wei Han}, Blow-up of solutions to semilinear wave equations with variable
coefficients and boundary, J. Math. Anal. Appl. ${\bf 374}$ $(2011),$ $585-601.$

\end{document}